\documentclass[pdflatex,sn-mathphys-num]{sn-jnl}

\usepackage{graphicx}%
\usepackage{multirow}%
\usepackage{amsmath,amssymb,amsfonts}%
\usepackage{amsthm}%
\usepackage{mathrsfs}%
\usepackage[title]{appendix}%
\usepackage{color,xcolor}
\usepackage{textcomp}%
\usepackage{manyfoot}%
\usepackage{booktabs}%
\usepackage{algorithm}%
\usepackage{algorithmicx}%
\usepackage{algpseudocode}%
\usepackage{listings}%
\usepackage{natbib}
\usepackage[table]{xcolor}
\usepackage{colortbl}
\usepackage{makecell}

\theoremstyle{thmstyleone}%
\theoremstyle{thmstyletwo}%

\theoremstyle{thmstylethree}%

\begin{document}

\title[Large Language Models for Operations Research: A Comprehensive Survey]{Large Language Models for Operations Research: A Comprehensive Survey}


\author[1]{\fnm{Xianchao} \sur{Xiu}}\email{xcxiu@shu.edu.cn}

\author[1]{\fnm{Jianhao} \sur{Li}}\email{lijianhao@shu.edu.cn}

\author*[2]{\fnm{Jun} \sur{Fan}}\email{junfan@hebut.edu.cn}

\author[3]{\fnm{Wanquan} \sur{Liu}}\email{liuwq63@mail.sysu.edu.cn}

\affil[1]{\orgdiv{School of Mechatronic Engineering and Automation}, \orgname{Shanghai University}, \orgaddress{\city{Shanghai}, \postcode{200444}, \country{China}}}

\affil*[2]{\orgdiv{Institute of Mathematics}, \orgname{Hebei University of Technology}, \orgaddress{\city{Tianjin}, \postcode{300401}, \country{China}}}

\affil[3]{\orgdiv{School of Intelligent Systems Engineering}, \orgname{Sun Yat-sen University}, \orgaddress{\city{Guangzhou}, \postcode{510275}, \country{China}}}


\abstract{Operations Research (OR) serves as a core decision-support methodology for complex systems, with significant applications across mathematics, management science, and computer science. Traditional approaches heavily rely on expert knowledge and often struggle to efficiently solve large-scale and multi-constraint problems. The rapid advancement of Large Language Models (LLMs) in recent years has offered a novel research paradigm to address these challenges. This paper presents a systematic survey of Large Language Models for Operations Research (LLM4OR). We begin by introducing the definition of OR problems and the fundamental principles of LLMs. We then focus on analyzing the roles of LLMs in OR, specifically covering such as model formulation, algorithm design, and solution verification. In addition, we discuss practical applications in representative scenarios and summarize benchmark datasets in this field. Finally, we outline the key challenges and provide perspectives on future research directions. A  collection of related literature is available at \url{https://github.com/xianchaoxiu/LLM4OR} for reference. 
}




\maketitle

\section{Introduction}\label{sec1}

Operations Research (OR) is an applied mathematical discipline that primarily focuses on using mathematical models and optimization approaches to solve complex problems. Over the past few decades, it has played a pivotal role in fields such as scheduling, transportation, robotics, biochemistry, networking, and finance \cite{hillier2005introduction}.

\begin{figure}[b]
	\centering
	\includegraphics[scale=0.34]{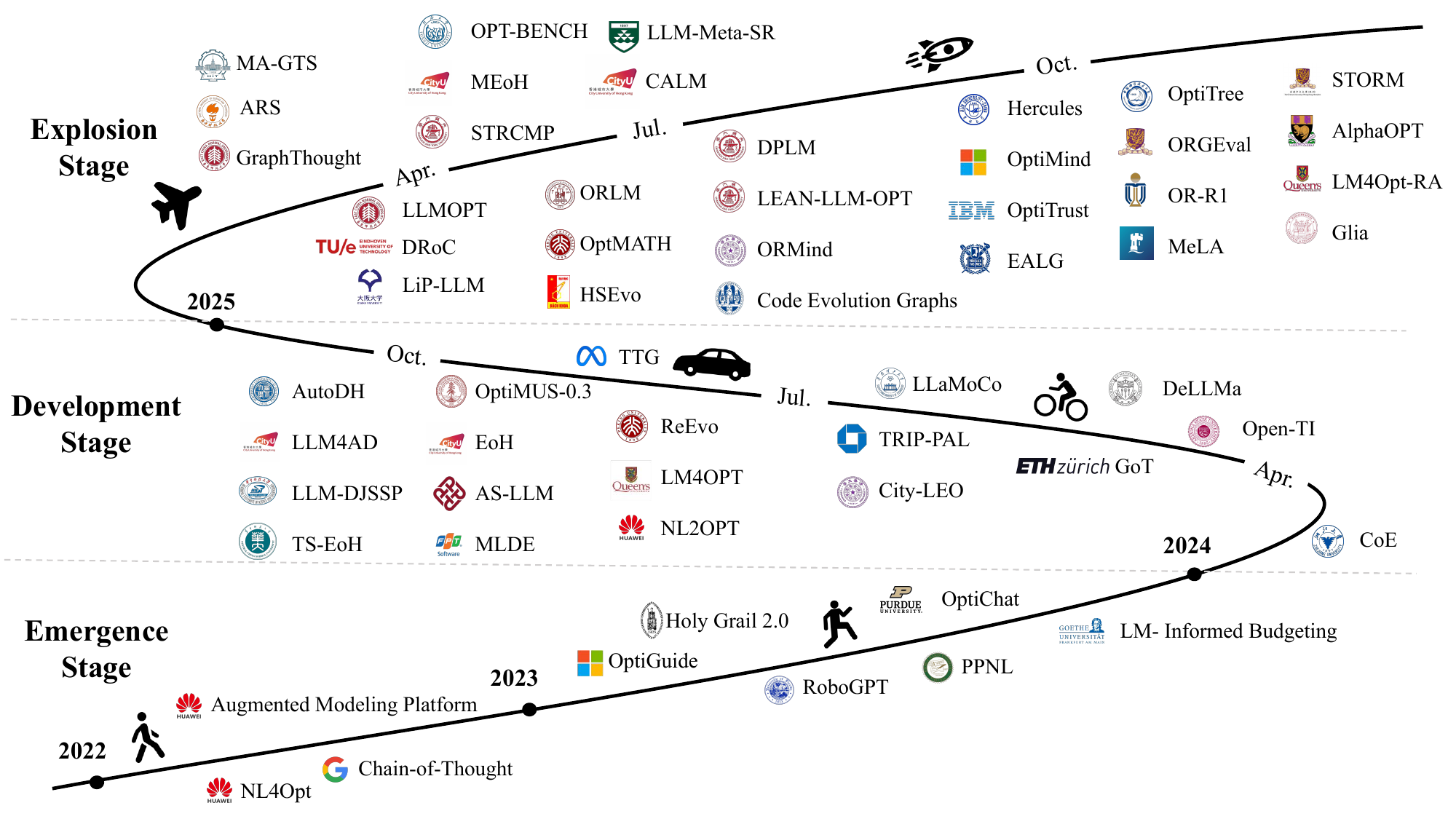} 
	\caption{Evolution timeline of LLM4OR (2022-2025).}
	\label{fig1}
\end{figure}

Large Language Models (LLMs), leveraging their massive parameter scales
and high-quality training data, have demonstrated strong generalization capabilities
across multiple tasks, marking an important step toward Artificial General Intelligence (AGI)~\cite{brown2020language, wang2022self}. 
In recent years, LLMs have shown broad potential in complex domains such as mathematical reasoning \cite{ahn2024large, li2026evosr}, code programming \cite{jiang2026survey} and algorithm design \cite{liu2025systematic}. They also achieve promising progress in data analysis \cite{sun2025lambda} and feature selection \cite{li2025llm4fs}, providing new perspectives for research in Artificial Intelligence for Science (AI4S). 
Against this background, the Large Language Models for Operations Research (LLM4OR) has attracted widespread attention from both academia and industry, giving rise to numerous representative studies as shown in Fig. \ref{fig1}. 

LLMs not only reduce the manpower required for problem solving, but also improve the efficiency of solutions, reshaping the research paradigm of OR problems to a certain extent. 
Fan et al. \cite{fan2025artificial} summarized the applications of Artificial Intelligence (AI) techniques such as neural networks and Reinforcement Learning (RL) in solving OR problems. 
Meanwhile, Huang et al. \cite{huang2024large} studied the mutually reinforcing interactions between LLMs and OR algorithms. 
Wu et al. \cite{wu2025evolutionary} explored the collaborative mechanisms between LLMs and evolutionary algorithms, demonstrating their roles in code generation and software engineering. 
In addition, Guo et al. \cite{tiande2025machine} provided an in-depth analysis of Machine Learning (ML) methods for solving combinatorial OR problems, while Forootani et al. \cite{forootani2025survey} reviewed the capabilities of LLMs in OR algorithms and examined their applications in mixed integer programming, linear quadratic control, and multi-agent frameworks. 
Besides, Xiao et al. \cite{xiao2025survey} focused on LLMs in the model formulation stage of OR problems and related benchmark datasets. 
Recently, Zhang et al. \cite{zhang2025systematic} reviewed how LLMs empower evolutionary OR algorithms from both modeling and solution perspectives. 
Wang et al. \cite{wang2025large} summarized LLM-based methods, applications, and challenges from the aspects of automated modeling, auxiliary optimization, and typical scenarios. 
Table \ref{tab1} provides a comparison of the current surveys on LLM4OR.

\begin{table}[t]
\centering
\caption{Comparison with other representative surveys.}
\label{tab1}
\renewcommand{\arraystretch}{1.2}
\footnotesize
\setlength{\tabcolsep}{4pt}
\begin{tabular*}{\textwidth}{@{\extracolsep\fill} c c c c c c}
\toprule
References & \makecell{Model \\ Formulation} & \makecell{Algorithm \\ Design} & \makecell{Solution \\ Validation} & \makecell{Scenario \\ Application} & \makecell{Benchmark \\ Datasets} \\
\midrule
\cite{fan2025artificial}        & $\checkmark$ & $\checkmark$ & $\times$     & $\times$     & $\times$     \\
\cite{huang2024large}           & $\times$     & $\checkmark$ & $\times$     & $\times$     & $\times$     \\
\cite{wu2025evolutionary}      & $\times$     & $\checkmark$ & $\times$     & $\checkmark$ & $\times$     \\
\cite{tiande2025machine}       & $\checkmark$ & $\checkmark$ & $\times$     & $\times$     & $\times$     \\
\cite{forootani2025survey}     & $\checkmark$ & $\checkmark$ & $\times$     & $\times$     & $\checkmark$ \\
\cite{xiao2025survey}          & $\checkmark$ & $\times$     & $\checkmark$ & $\times$     & $\checkmark$ \\
\cite{zhang2025systematic}     & $\checkmark$ & $\checkmark$ & $\times$     & $\checkmark$ & $\times$     \\
\cite{wang2025large}           & $\checkmark$ & $\checkmark$ & $\times$     & $\checkmark$ & $\checkmark$ \\ \midrule
This Work                     & $\checkmark$ & $\checkmark$ & $\checkmark$ & $\checkmark$ & $\checkmark$ \\
\bottomrule
\end{tabular*}
\end{table}

Given the rapid development of LLM4OR, there is an urgent need to systematically summarize the research progress in this field. 
To this end, we conduct a survey of the literature published up to May 2026 and complete this work, whose main contributions are summarized as follows:
\begin{itemize}
\item We provide a comprehensive survey of the current status and development trends in LLM4OR, summarizing nearly 180 highly relevant publications and establishing an evolving public repository.
\item To the best of our knowledge, this is the first work to comprehensively cover key steps such as model construction, algorithm design, and solution verification, thereby aligning itself more closely with the inherent logic of OR.
\item We explore the real-world application scenarios and various benchmark datasets in LLM4OR, as well as existing challenges, providing a valuable reference for both academia and industry.
\end{itemize}

The structure of the paper is as follows.
Section \ref{sec2} briefly introduces the background of OR problems as well as LLMs.
Section \ref{sec3}, Section \ref{sec4}, and Section \ref{sec5} specifically show the roles of LLMs in model formulation, algorithm design, and solution validation,  respectively. Section \ref{sec6} presents typical application scenarios. Section \ref{sec7} elaborates on relevant datasets. 
Section \ref{sec8} provides conclusions and perspectives for future research.

\section{Background}\label{sec2}

Based on the types of decision variables, OR problems can be classified into two major categories: combinatorial optimization and continuous optimization. 
This scetion introduces the basic definitions, characteristics, and algorithms of these two types, and then reviews the fundamental principles of LLMs.

\subsection{Operations Research Problems}
\subsubsection[Combinatorial Optimization]{Combinatorial Optimization}

Combinatorial optimization problems involve seeking optimal solutions within a discrete feasible domain. 
Typical examples include the Traveling Salesman Problems (TSP), Graph Coloring Problems (GCP), and Scheduling Problems (SP). 
The mathematical model of a combinatorial optimization problem can be expressed as
\begin{equation}
\begin{aligned}
\min_{\mathbf{x} \in \mathbb{Z}^n} \quad & f(\mathbf{x}) \\
\text{s.t.} \quad 
& g_i(\mathbf{x}) \le 0, \quad i = 1, 2, \dots, m, \\
& h_j(\mathbf{x}) = 0, \quad j = 1, 2, \dots, p,
\end{aligned}
\label{eq1}
\end{equation}
where $\mathbf{x} \in \mathbb{Z}^n$ denotes the discrete decision variable, $f(\mathbf{x})$ is the objective function, $g_i(\mathbf{x})$ and $h_j(\mathbf{x})$ represent the inequality and equality constraints, respectively.

Due to the discrete nature of $\mathbf{x}$, the number of feasible solutions in combinatorial optimization problems often grows exponentially with the problem dimension, and a large number of local optimal solutions may exist. 
Therefore, in terms of computational complexity, many classical models belong to NP-hard problems. 
Combinatorial optimization algorithms can generally be categorized as exact algorithms  and heuristic algorithms. 
Exact algorithms guarantee convergence to the global optimum, such as exhaustive search methods, cutting‑plane methods, and branch‑and‑bound methods. 
However, they often suffer from high computational complexity for large‑scale problems. 
In contrast, heuristic algorithms aim to find high‑quality feasible solutions within a reasonable computational time, such as greedy methods, simulated annealing methods, and genetic methods. 
Although these algorithms generally cannot obtain global optimality, they demonstrate strong practicality in large‑scale and complex engineering scenarios.

\subsubsection[Continuous Optimization]{Continuous Optimization}

Unlike combinatorial optimization problems, continuous optimization problems aim to find optimal solutions within a continuous feasible domain, and they have important applications in multiple fields such as ML, control systems, and signal processing. 
The mathematical model of a continuous optimization problem can be characterized by
\begin{equation}
\begin{aligned}
\min_{\mathbf{x} \in \mathbb{R}^n} \quad & f(\mathbf{x}) \\
\text{s.t.} \quad 
& g_i(\mathbf{x}) \le 0, \quad i = 1, 2, \dots, m, \\
& h_j(\mathbf{x}) = 0, \quad j = 1, 2, \dots, p,
\end{aligned}
\label{eq2}
\end{equation}
where $\mathbf{x} \in \mathbb{R}^n$ represents the continuous decision variable.

Based on the properties of the objective and constraint functions, the continuous optimization problem \eqref{eq2} can be classified into convex optimization and nonconvex optimization. 
When $f$ and $\{ g_i \}$ are convex functions and $\{ h_j \}$ are affine functions, problem \eqref{eq2} is a convex optimization problem, for which algorithms can be designed to converge to the global optimal solution under appropriate conditions. 
In contrast, for nonconvex optimization, it is generally difficult to guarantee obtaining the global optimal solution, and algorithms typically converge to stationary points or local optima. Common numerical algorithms  for continuous optimization include gradient descent methods, Newton's methods, trust-region methods, and augmented Lagrangian methods. Interested readers may refer to \cite{nocedal2006numerical}.

\subsection[Large Language Models]{Large Language Models}
\subsubsection[Basic Principles]{Basic Principles}

The core architecture of LLMs is the Transformer~\cite{vaswani2017attention}. 
Through the self-attention mechanism, it can capture long-range dependencies and assign weights to the correlations among tokens, thereby efficiently processing large-scale datasets. 
In the Transformer architecture, the encoder and decoder play different roles in representation learning and sequence generation, respectively. The encoder maps the input sequence into a set of context-dependent high-dimensional representations, typically consisting of multiple stacked self-attention sublayers and feedforward neural networks. 
The decoder generates the target sequence based on the given context. Its structure is similar to that of the encoder, but its self-attention mechanism usually employs a causal mask to ensure that the current token depends only on previously generated tokens. 
In the standard Transformer architecture, the decoder interacts with encoder outputs through cross-attention mechanisms to establish mappings between input and output sequences. Nevertheless, most mainstream text generation models adopt a simpler decoder-only structure. 
Representative LLMs including the Generative Pre-trained Transformer (GPT) series and Large Language Model Meta AI (LLaMA) series are fully constructed relying on autoregressive decoders.

\subsubsection[Training Paradigm]{Training Paradigm}

The training of LLMs can be broadly divided into two key stages: pre-training~\cite{min2023recent} and fine-tuning~\cite{ding2023parameter}. 
Pre-training is conducted on large-scale unlabeled text corpora, enabling LLMs to learn linguistic patterns and generate logically coherent text. 
Representative pre-training approaches include autoregressive models and masked language models. 
The GPT series adopts an autoregressive paradigm, where each generated token depends on the preceding tokens. 
In contrast, models such as Bidirectional Encoder Representations from Transformers (BERT) employ the masked language modeling approach, where a subset of tokens in the input text is randomly masked, and the model is trained to predict the masked tokens.

Fine-tuning further retrains the pre-trained LLMs on specific tasks or datasets. 
Through fine-tuning, LLMs can adjust their internal parameters according to task requirements, thereby effectively transforming general-purpose models into task-specific models aligned with human preferences. 
The fine-tuning process typically requires a certain scale of labeled datasets and appropriate loss functions depending on the task, such as cross-entropy loss for classification tasks and negative log-likelihood loss for generation tasks. 
fine-tuning encompasses several strategies, including Supervised Fine-Tuning (SFT)~\cite{dong2024abilities}, Reinforcement Learning from Human Feedback (RLHF)~\cite{kaufmann2024survey}, and Parameter-Efficient Fine-Tuning (PEFT)~\cite{han2024parameter}, which will be discussed in subsequent sections. 
SFT trains LLMs using high-quality human-annotated data, enabling them to generate appropriate outputs based on given instructions. 
RLHF further incorporates human preferences by introducing reward models and RL algorithms to align model outputs with desired values. 
To reduce the computational and storage costs of full-parameter fine-tuning, PEFT have been proposed. 
For example, Low-Rank Adaptation (LoRA) ~\cite{hu2022lora} injects parameters through low-rank adaptation matrices, while prefix-tuning~\cite{li2021prefix} achieves efficient adaptation by introducing trainable prefixes. 
These fine-tuning strategies are not strictly sequential stages, but can be applied independently or in combination.

\begin{figure}[t]
	\centering
	\includegraphics[scale=0.43]{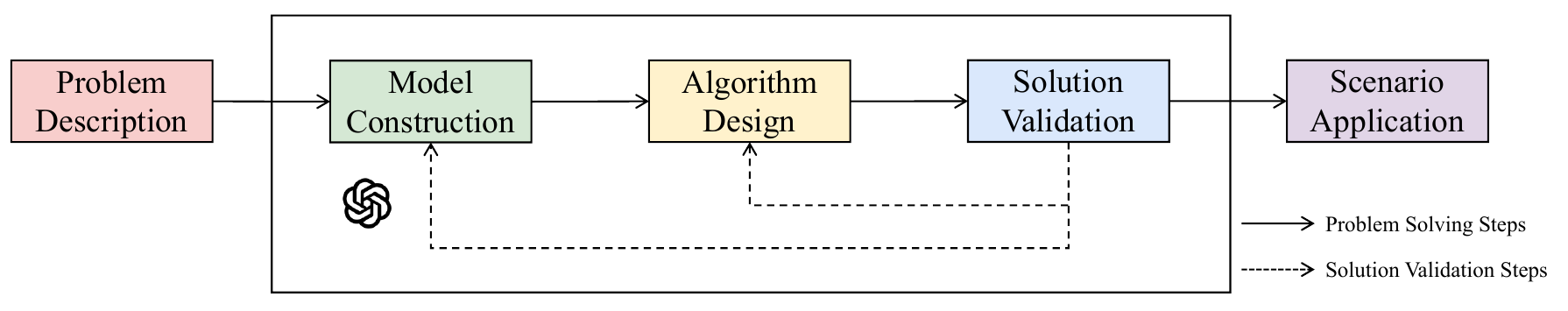} 
	\caption{Solution workflow for OR problems.}
	\label{fig3}
\end{figure}

\subsubsection[Solution Process]{Solution Process}

In addition to problem description and scenario application, the solution process of OR problems can be summarized into three steps:
\begin{itemize}
    \item \textit{Model Formulation}: This step is to describe the text problem, including  objective functions and constraints, then to construct an accurate mathematical model. 
    \item \textit{Algorithm Design}: Once the model is established, appropriate optimization algorithms can be applied according to its type and complexity.
    \item \textit{Solution Validation}: The proposed model and algorithms are evaluated to determine their effectiveness, and, if necessary, further refine them.
\end{itemize}
As illustrated in Fig.~\ref{fig3}, LLMs can play a key role in each stage of the OR problem solving process, thereby enhancing the overall level of automation.

\section{Model Formulation}\label{sec3}

Model formulation is a fundamental step in OR, aiming to formalize problems described in natural languages into mathematical models. 
The quality of the constructed model largely influences subsequent algorithm selection and solution efficiency, and ultimately determines the accuracy. 
Traditional approaches rely heavily on expert knowledge, whereas the strong mathematical reasoning capabilities of LLMs provide new perspectives for this process~\cite{azerbayev2024llemma,yang2024formal}. 
LLMs are capable of extracting problem background, objective functions, and constraints from natural language descriptions, and transforming them into mathematical formulations.

Ramamonjison et al.~\cite{ramamonjison2022augmenting} developed an intelligent modeling platform based on the BERT language model, enabling users to interact with AI systems to efficiently construct models for OR problems. 
The team also organized the NL4Opt competition, which aims to leverage LLMs for the modeling process. 
It consists of two challenging subtasks: (1) identifying objective functions, constraints, and decision variables, (2) constructing mathematical models based on the extracted key elements~\cite{ramamonjison2023nl4opt}. In addition, Ahmed et al.~\cite{ahmed2024lm4opt} put forward LM4OPT to evaluate the mathematical modeling performance of mainstream LLMs like GPT-3.5, GPT-4 and Llama-2-7B.
For robust and adaptive optimization research, Bertsimas et al.~\cite{bertsimas2024robust} verified that ChatGPT is competent to complete modeling and solving work for these problems.
In addition, Zhai et al.~\cite{zhai2025equivamap}  constructed EquivaMap to automatically find equivalent formulations of optimization tasks. Overall, the role of LLMs in the model formulation of OR problems has grown increasingly prominent. 
Two main research paradigms have gradually emerged in this direction: prompt-based methods and learning-based methods.

\begin{figure}[t]
	\centering
	\includegraphics[scale=0.38]{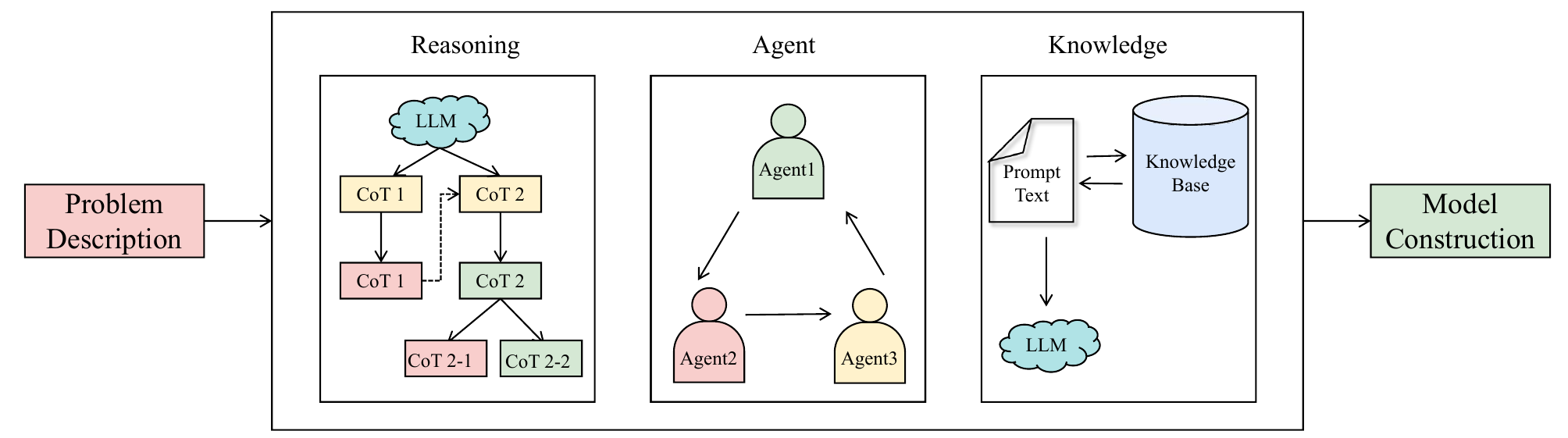} 
	\caption{Illustration of three frameworks for prompt-based methods.}
	\label{fig5}
\end{figure}

\subsection{Prompt-Based Methods}

Prompt-based methods guide LLMs to perform tasks more effectively without relying on large-scale parameter tuning, by designing specialized prompts and leveraging techniques such as multi-agent frameworks and Retrieval-Augmented Generation (RAG). 
This approach represents a typical research paradigm in natural language processing. 
It fully exploits the general knowledge acquired by LLMs during pre-training, allowing them to demonstrate strong zero-shot or few-shot capabilities based solely on the provided prompts~\cite{yao2023react}. 
Prompt-based methods can be categorized into three frameworks: reasoning, agent, and knowledge, as illustrated in Fig.~\ref{fig5}.

\subsubsection[Reasoning Framework]{Reasoning Framework}

The reasoning framework aims to enable LLMs to perform effective reasoning, thereby proposing more appropriate solutions, including Chain-of-Thought (CoT)~\cite{wei2022chain}, Tree-of-Thought (ToT)~\cite{yao2023tree}, and Graph-of-Thought (GoT)~\cite{besta2024graph}. 
Although these methods were originally designed for general tasks, they have also been successfully applied to model formulation~\cite{xiao2024chain}. 
Inspired by the NL4Opt competition, Tsouros et al.~\cite{tsouros2023holy} developed a complete framework spanning from mathematical modeling to solvers. 
For higher-quality results, a recursive dynamic temperature strategy combined with CoT was proposed to acquire better feasible solutions~\cite{wang2024leveraging}. 
In view of hierarchical analysis, a modeling tree sorted by problem hierarchy and complexity was built, so that OR problems can be matched with the most similar problem templates accordingly~\cite{liuoptitree}.


\subsubsection[Agent Framework]{Agent Framework}

\begin{figure}[t]
	\centering
	\includegraphics[scale=0.43]{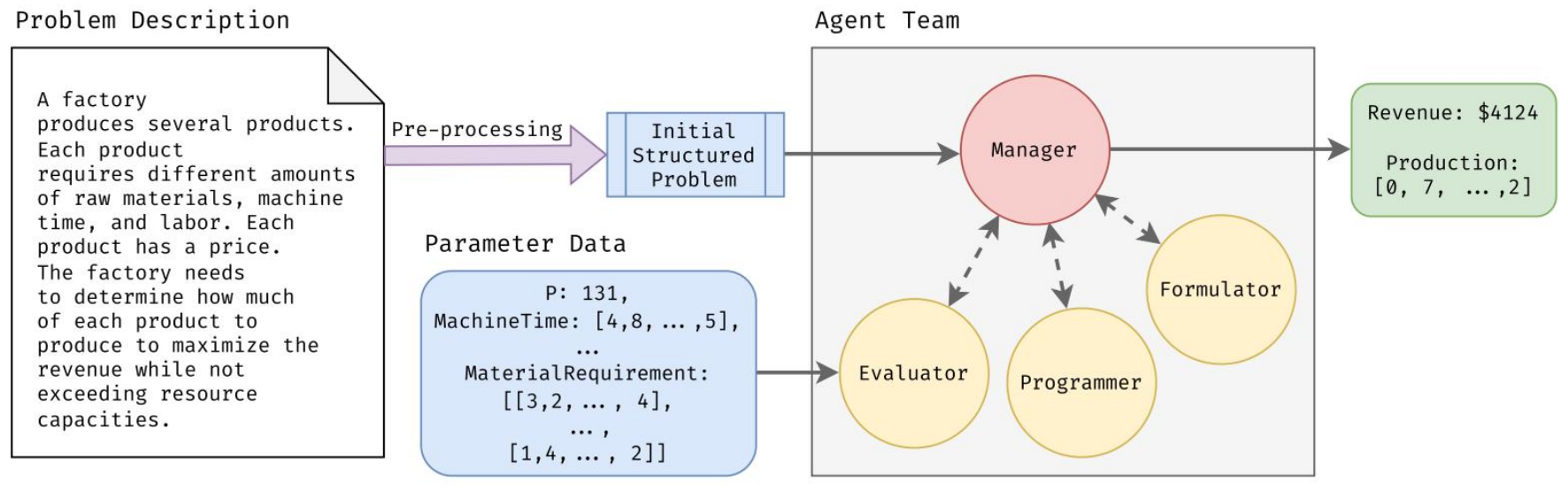} 
	\caption{Workflow of OptiMUS \cite{ahmaditeshnizi2024optimus0.2}.}
	\label{fig6}
\end{figure}

In the agent framework, LLMs are assigned to multiple agents (e.g., formula generator, code generator, evaluator), allowing each agent to work independently while coordinating with others. 
To efficiently handle long descriptions and complex data while mitigating prompt bloat, a modular system named OptiMUS 0.2 was designed  \cite{ahmaditeshnizi2024optimus0.2}. 
By implementing cross-validation among agents, this approach significantly elevates the mathematical modeling capabilities of LLMs, as illustrated in Fig.~\ref{fig6}. 
The effectiveness of such systems is further underscored by LLMFP \cite{hao2025planning} and validated across nine distinct tasks using two different LLMs. 
Building on these specialized structures, various auxiliary mechanisms have been introduced to further refine modeling performance. For instance, Xiao et al.~\cite{xiao2024chain} proposed a method termed CoE, which incorporates forward reasoning and backward reflection to coordinate inter-module tasks. 
Mostajabdaveh et al.~\cite{mostajabdaveh2024optimization} leveraged collaborative and competitive strategies among agents to validate results, thereby enhancing robustness in solving complex problems. 
In terms of error correction, OptimAI~\cite{thind2025optimai} was developed to dynamically switch among alternative solutions during the debugging phase. Besides, Monte Carlo Tree Search (MCTS) was integrated into an automated formulation framework, which employs pruning techniques to eliminate trivial equivalent formulas~\cite{astorga2025autoformulation}. 
To address high computational latency and limited synergy inherent in sequential coordination, Berto et al.~\cite{berto2025parco} designed a parallel autoregressive paradigm named PARCO. This paradigm integrates a Transformer communication layer, a multi-pointer mechanism, and priority-based conflict, demonstrating strong performance in complex operational research tasks like path planning, pickup-and-delivery, and job shop scheduling.

\subsubsection[Knowledge Framework]{Knowledge Framework}

RAG techniques provide LLMs with additional support from external knowledge bases, significantly enhancing the performance of OR problem solving~\cite{jiang2025large}. 
Building upon earlier collaborative designs, an integration of RAG was introduced in OptiMUS 0.3, effectively enhancing the accuracy of its predecessor~\cite{ahmaditeshnizi2024optimus0.3}. 
To tackle complex mathematical formulations, Jiang et al.~\cite{jiang2025droc} proposed a constraint-decomposition retrieval method termed DROC, which is capable of breaking down intricate constraints to reduce overall modeling complexity. 
Focusing on practical deployment, a design combining local LLMs with domain-specific knowledge bases was developed and successfully applied to an aircraft skin manufacturing case~\cite{peng2025automatic}. 
Leveraging few-shot prompting alongside collaborative execution mechanisms, Liang et al.~\cite{liang2026llm} incorporated RAG into LEAN-LLM-OPT, which enables efficient solutions for complex OR problems involving long textual descriptions and external data inputs.


\subsection[Learning-Based Methods]{Learning-Based Methods}

\begin{figure}[t]
	\centering
	\includegraphics[scale=0.38]{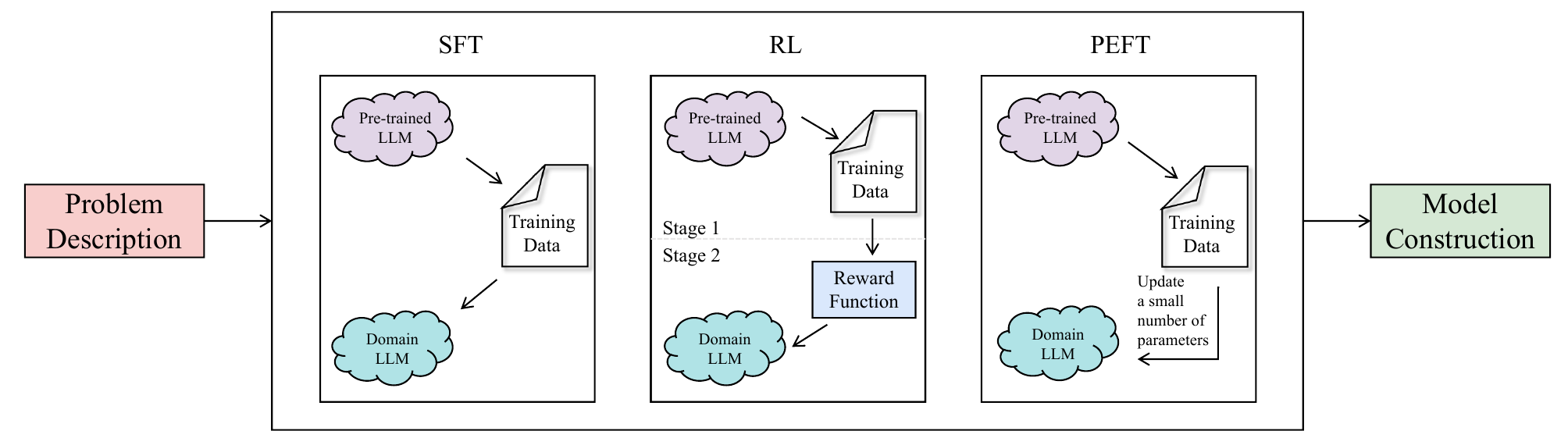} 
	\caption{Illustration of three frameworks for learning-based methods.}

	\label{fig7}
\end{figure}

In contrast to prompt-based methods, learning-based methods further adapt pre-trained LLMs to specific tasks and domains by incorporating strategies such as SFT, RL, and PEFT, which are illustrated in Fig.~\ref{fig7}.

\subsubsection[Supervised Fine-Tuning]{Supervised Fine-Tuning}

Under the SFT paradigm, researchers have conducted extensive explorations on data construction and training strategies. 
Amarasinghe et al.~\cite{amarasinghe2023ai} proposed AI Copilot, which leverages an LLM fine-tuned method via SFT for modeling by designing nine submodules and applying prompt engineering strategies to mitigate the token limitation of LLMs. 
Based on the NL4Opt dataset, Li et al.~\cite{li2023synthesizing} expanded the range of problem descriptions and constraint types, thereby significantly improving modeling accuracy by fine-tuning ChatGPT and Google Bard on this augmented dataset. 
The effectiveness of fine-tuned models was also demonstrated on the TSP, where a self-ensemble method was used to further refine solutions~\cite{masoud2024exploring}. 
Beyond directly expanding datasets, Yang et al.~\cite{yang2025optibench} proposed ReSocratic, a reverse data synthesis method from mathematical models to problem descriptions, primarily using the generated dataset to train multiple open-source models. 
In a similar vein, Huang et al.~\cite{huang2025orlm} introduced a data synthesis method named OR-INSTRUCT, constructing datasets via two distinct strategies. 
The first is expansion, which increases data coverage by extending problem scenarios and types, while the second is augmentation, which enhances data diversity by rewriting objective functions and constraints, restating problems, and introducing multiple modeling techniques, ultimately resulting in improved performance across several models. 
Additionally, Ma et al.~\cite{ma2026llamoco} introduced LLaMoCo, which incorporates Contrastive Learning (CL) as a warm-up phase to improve training efficiency and stability during fine-tuning. 
Besides, Jiang et al.~\cite{jiang2026large} proposed a two-stage strategy combining SFT and Feasibility-and-Optimality-Aware Reinforcement Learning (FOARL) to guide LLMs and improve solution quality.

\subsubsection[Reinforcement Learning]{Reinforcement Learning}

RL strategies are introduced to enhance model robustness. 
To address hallucination issues in LLMs, Jiang et al.~\cite{jiang2025llmopt} incorporated Kahneman-Tversky Optimization (KTO)~\cite{ethayarajh2024kto} along with self-correction mechanisms, and proposed LLMOPT, which has been validated across six real-world datasets spanning 20 domains, including healthcare, environment, and manufacturing. 
Similarly, Zhou et al.~\cite{zhou2025auto} proposed DPLM, which employs a two-stage training procedure: the first stage uses SFT, and the second stage uses either Direct Preference Optimization (DPO)~\cite{rafailov2023direct} for offline training or Group Relative Policy Optimization (GRPO)~\cite{shao2024deepseekmath} for online training, achieving strong performance in modeling dynamic programming tasks. 
Likewise, Ding et al.~\cite{ding2026or} combined SFT with Test-time Group Relative Policy Optimization (TGRPO) to develop OR-R1, enabling LLMs to effectively leverage both scarce labeled data and abundant unlabeled data during training.

\subsubsection[Parameter-Efficient Fine-Tuning]{Parameter-Efficient Fine-Tuning}

Considering that SFT requires updating all model parameters, which often incurs high time and computational costs, PEFT has attracted significant attention. 
Among these methods, LoRA is currently the most widely adopted approach. 
To validate and filter synthetic data, OptMATH~\cite{lu2025optmath} was equipped with a suppression-sampling mechanism, sustaining stable performance among LoRA-fine-tuned models ranging from 0.5B to 32B parameters. 
Zhang et al.~\cite{zhang2024solving} introduced OptLLM, which leverages external solvers to assist decision-making and supports multi-turn interactive modeling, significantly improving the modeling accuracy of the Qwen model after LoRA. 
Wu et al.~\cite{wu2025evostep} designed EVO-STEP-INDUCT, which guides LLMs to generate high-quality, diverse data through dual strategies of complexity evolution and range evolution, demonstrating the effectiveness of combining problem-structure validation with stepwise improvement in LoRA experiments on Llama-3-8B and Mistral-7B.\\

\textit{Prompt-based methods fully leverage the reasoning capabilities of LLMs. Through the reasoning framework, agent framework, and knowledge framework, LLMs can perform model formulation for OR problems with little or no additional training, offering advantages such as low implementation cost and flexible deployment, and enabling rapid transfer to different types of OR tasks. 
In contrast, learning-based methods, relying on high-quality data synthesis strategies and various training mechanisms, drive LLMs to learn modeling patterns from data through SFT, RL, and PEFT, thereby improving both accuracy and stability.}

\section[Algorithm Design]{Algorithm Design}\label{sec4}

\begin{figure}[t]
	\centering
	\includegraphics[scale=0.38]{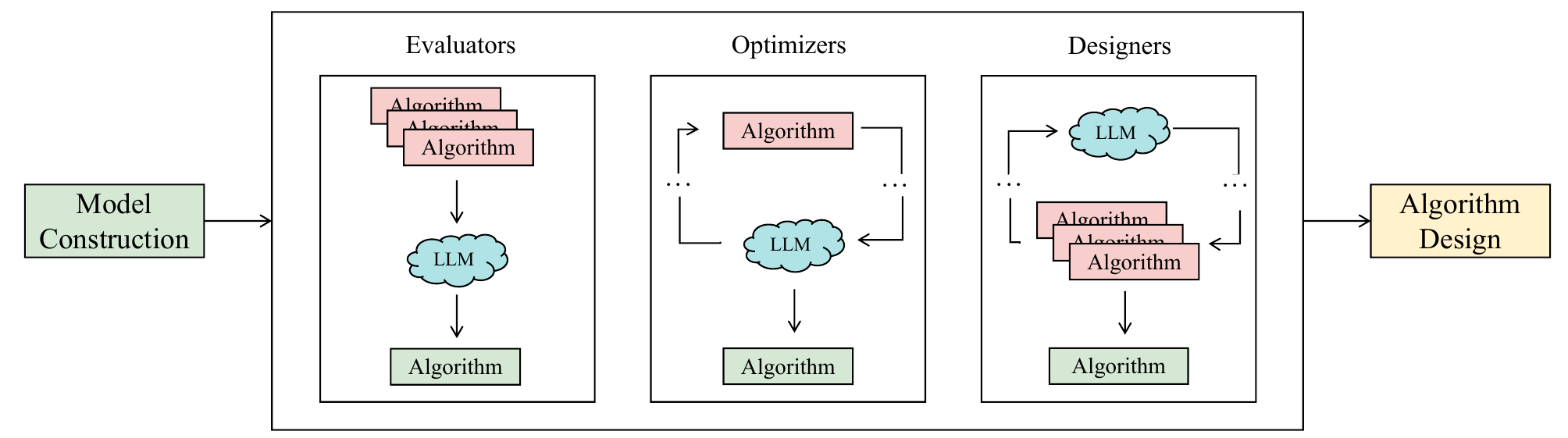} 
    \caption{Illustration of three roles in algorithm design.}

	\label{fig9}
\end{figure}

Algorithm design is a critical component of OR, directly affecting both the efficiency and accuracy of problem solving~\cite{yao2024evolve, lou2026llm}. 
Researchers have begun incorporating LLMs into the automated algorithm design process, aiming to improve both the efficiency and quality of algorithm generation. 
This trend shifts OR algorithms from being experience-driven toward LLM-driven generation, not only lowering the barrier for algorithm development but also providing new perspectives for constructing general-purpose algorithmic frameworks. LLM-driven algorithms have demonstrated excellent performance across a variety of classical OR problems, sometimes surpassing expert-designed algorithms, while exhibiting strong generalization capabilities~\cite{lehman2023evolution, xu2026synergistic}. 
Based on the roles that LLMs play in algorithm design, they can be categorized into three types: evaluators, optimizers, and designers, as illustrated in Fig.~\ref{fig9}. 

\subsection[Combinatorial Optimization]{Combinatorial Optimization}

In combinatorial optimization, LLMs can improve search strategies and also integrate with traditional optimization methods~\cite{liu2025fitness, da2026large}.

\subsubsection[Evaluators]{Evaluators}

Evaluators refer to LLMs serving as auxiliary modules in OR, used for tasks such as algorithm evaluation and selection. 
Nie et al.~\cite{nie2023importance} designed an LLM-based solver that synthesizes directional feedback from historical optimization trajectories (similar to first-order derivative information in traditional optimization methods) to achieve reliable improvements during iterative processes. 
Wu et al.~\cite{wu2024large} proposed AS-LLM by making LLMs to capture the structural and semantic characteristics of algorithms. 
Similarly, Li et al.~\cite{li2025strcmp} introduced STRCMP by combining graph neural networks with LLMs, enabling LLMs to extract structural embeddings from combinatorial optimization problem instances to effectively identify high-performing algorithms.
Additionally, Wu et al.~\cite{wu2025efficient} proposed the few-shot performance prediction prompting framework called Hercules, which analyzes the semantic similarity between LLM-generated algorithms and existing algorithms, enabling rapid evaluation of heuristic algorithm performance.

\subsubsection[Optimizers]{Optimizers}

Optimizers refer to LLMs that actively search, mutate, or recombine existing algorithms. 
Mao et al.~\cite{mao2024identify} guided LLMs to evolve manually designed node scoring functions, and the generated functions demonstrated strong adaptability in decision tasks for identifying key nodes in networks. 
Yu et al.~\cite{yu2024autornet} combined LLMs with evolutionary algorithms to develop AutoRNet for the design of robust networks. 
Some studies have incorporated fine-tuning and other techniques to enhance the capability of LLMs in improving algorithms. 
For example, Zhang et al.~\cite{zhang2024gcoder} adopted a two-stage training approach via combining SFT with Reinforcement Learning from Compiler Feedback (RLCF). 
Surina et al.~\cite{surina2025algorithm} integrated RL with evolutionary search to drive LLM-based algorithm improvements. 
Sartori et al.~\cite{sartori2025combinatorial,sartori2025improving}  demonstrated the potential of LLMs to enhance existing algorithms by validating across ten classical methods including metaheuristics, deterministic, and exact algorithms.

\begin{figure}[t]
	\centering
	\includegraphics[scale=0.25]{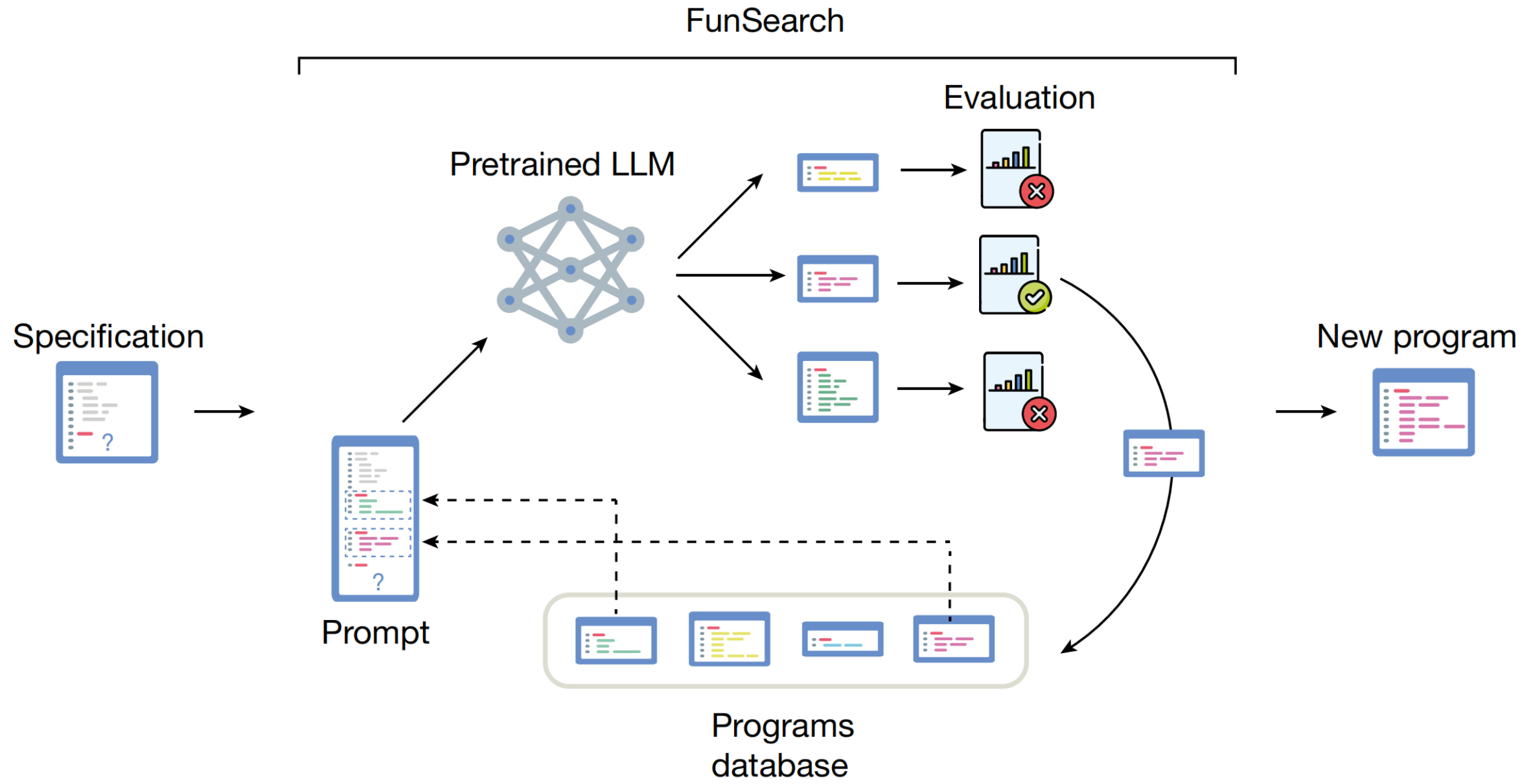} 
	\caption{Workfolw of FunSearch  \cite{romera2024mathematical}.}

	\label{FunSearch}
\end{figure}

\subsubsection[Designers]{Designers}

Designers refer to LLMs that directly design algorithms. 
Romera et al.~\cite{romera2024mathematical} combined pre-trained language models with a system evaluator to develop FunSearch, which performs searches over function spaces and achieves strong performance on the Cap Set Problems (CSP) and Online Bin Packing Problems (OBPP), as illustrated in Fig.~\ref{FunSearch}. 
Building on FunSearch, Chen et al.~\cite{chen2024qube} incorporated a quality-uncertainty balancing evolutionary strategy to design heuristic algorithms for NP-hard problems. 
Additionally, Liu et al.~\cite{liu2024evolution} proposed a framework called EOH through combining LLMs with evolutionary computation, aimed at enabling the model to autonomously design algorithms and their code. 
This framework has been extended to multiple domains. For example, Bomer et al.~\cite{bomer2025leveraging} developed CEOH for the Unit-Load Pre-Marshalling Problems (ULPMP), and Yao et al.~\cite{yao2025multi} proposed MEOH for the TSP and OBPP. 
To efficiently explore heuristic spaces, ReEvo~\cite{ye2024reevo} was constructed with embedded evolutionary search and LLM reflection mechanisms. 
For algorithm optimization, another improved heuristic approach is HSEvo~\cite{dat2025hsevo}, which maintains a trade-off between stability and population diversity, outperforming existing methods including EOH, FunSearch, and ReEvo. 
To provide a unified interface for algorithm design with LLMs, Liu et al.~\cite{liu2024llm4ad} proposed LLM4AD, capable of generating solving algorithms for multiple types of tasks. 
Yu et al.~\cite{yu2024deep} conducted in-depth analyses of three key components in algorithm design, i.e., individual representation, mutation operators, and fitness evaluation, while also combining LLMs with evolutionary algorithms to design heuristic algorithms. 
Sun et al.~\cite{sun2024autosat} demonstrated that LLMs can directly search for new heuristics for Satisfiability Problems (SATP). 
Shi et al.~\cite{shi2025generalizable} proposed the meta-heuristic optimization strategy called MOH, which automatically constructs diverse algorithms via meta-learning, showing strong performance across multiple combinatorial optimization problems. 
Furthermore, Huang et al.~\cite{huang2025calm} combined language-guided generation with numerical evaluation, and fine-tuned the LLMs with RL to enhance the quality of generated heuristic algorithms.

\subsection[Continuous Optimization]{Continuous Optimization}

For continuous optimization, LLMs can also be categorized into evaluators, optimizers, and designers. 
In the role of evaluators, Hao et al.~\cite{hao2024large} formulated model-assisted selection as a classification or regression task, where LLMs assess new solutions using historical data, leading to LAEA. 
Some research has focused on the role of optimizers. 
Wang et al.~\cite{wang2024large} proposed CMOEA-LLM, a constraint-based multi-objective optimization strategy that uses prompt engineering to guide LLMs in evaluating candidate quality and generating improved solutions. 
Guo et al.~\cite{guo2024llm} introduced an alternating gradient-based optimizer and LLM-based optimizer strategy to address complex nonconvex optimization problems.

Research surrounding the role of designers is rapidly advancing. 
Van et al.~\cite{van2024llamea} put forward a novel approach termed LLaMEA, capable of autonomously constructing algorithms to tackle Black-Box Optimization Problems (BBOP). 
On the basis of this approach, HyperParameter Optimization (HPO) is integrated into LLaMEA, yielding the improved variant LLaMEA-HPO~\cite{van2024loop} that has been deployed for BBOP, TSP and OBPP scenarios. 
Against the background of evolutionary optimization, EvoLLM~\cite{lange2024large} is developed to rank discretized individuals within populations in ascending order, and its outstanding effectiveness is verified across a BBOP benchmark consisting of 24 continuous optimization tasks as well as lightweight neural evolution experiments. 
Zhong et al.~\cite{zhong2024leveraging} utilized ChatGPT-3.5 to create a metaheuristic algorithm inspired by animal group behaviors, validating its effectiveness on six continuous optimization engineering problems. 
Brahmachary et al.~\cite{brahmachary2025large} proposed LEO and reliably validated on several industrial engineering problems, including supersonic nozzle shape optimization, heat transfer, and wind field layout optimization. 
Additionally, Liu et al.~\cite{liu2025large} introduced MOEA and D-LMO, using LLMs as black-box search operators in multi-objective evolutionary algorithms, generating candidate solutions, and further designing white-box operators to explain and approximate the behavior of LLMs.

After the algorithm structure is determined, parameter configuration is equally critical in the OR process, as it directly affects the efficiency of solution finding. 
With the ability to understand problem parameters, constraints, and objective functions, LLMs can automatically generate reasonable parameter configurations and continuously improve these strategies through feedback loops~\cite{ustyugov2024different, martinek2024large}. 
Studies have shown that even under conditions of limited search budgets (e.g., computational resources, time) or data scarcity, they can still achieve performance comparable to or even better than traditional parameter tuning methods~\cite{zhang2023using, lawless2025llms}. \\

\textit{Overall, by taking on roles such as evaluators, optimizers, and designers during algorithm design,  LLMs help alleviate the reliance on human expertise in OR algorithms, exhibiting strong generalization abilities in complex combinatorial and continuous optimization tasks.}

\begin{figure}[t]
	\centering
	\includegraphics[scale=0.38]{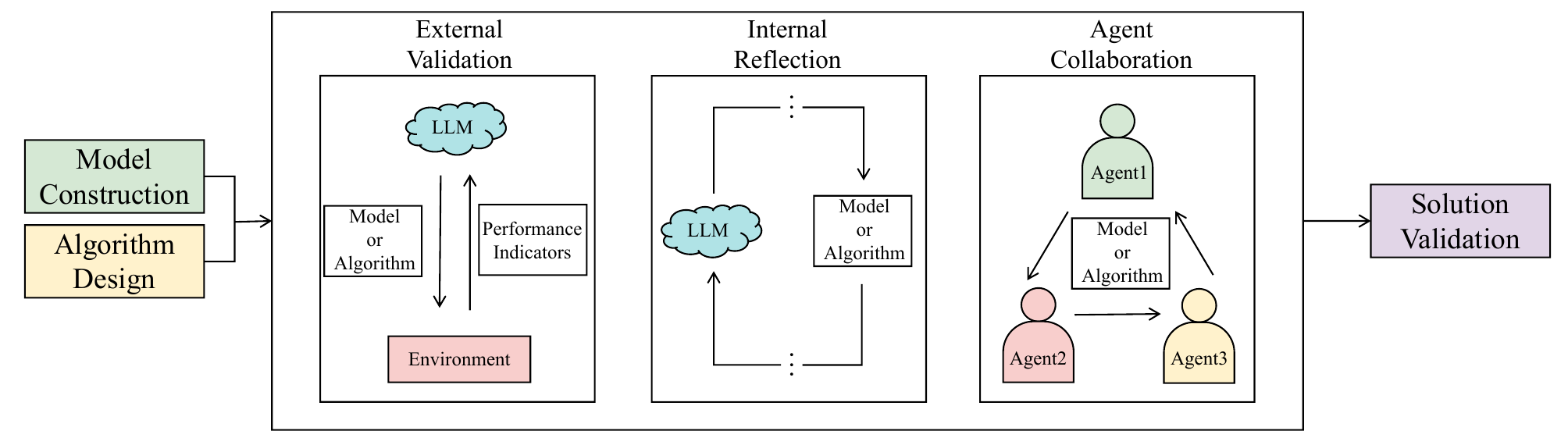} 
	\caption{Illustration of three solution validation frameworks.}

	\label{fig10}
\end{figure}

\section[Solution Validation]{Solution Validation}\label{sec5}

Across the complete workflow for OR problem solving, solution validation is a critical step to ensure that the model formulation and algorithm design meet practical requirements. 
Compared to traditional methods, while the models and algorithms output by LLMs are highly automated, they may suffer from issues such as constraint omissions or logical inconsistencies. 
Therefore, establishing an efficient solution validation mechanism is crucial for ensuring the reliability of LLM4OR. 
Current research mainly includes three frameworks: external validation, internal reflection, and agent collaboration, as shown in Fig.~\ref{fig10}.

\subsection[External Validation]{External Validation}

During the external validation phase, researchers typically employ the interaction feedback between LLMs and optimization solvers or RL reward signals as a validation mechanism to assess the generated solutions. 
Based on pre-trained language models, Almonacid et al.~\cite{almonacid2023towards} proposed an automatic modeling framework, utilizing two GPT systems for model generation and debugging, and interacting with the MINIZINC for error correction. 
Chen et al.~\cite{chen2024diagnosing} focused on the feasibility diagnosis of mixed-integer linear programming models, leveraging interactions between LLMs and the Gurobi to identify irreducibly infeasible problem subsets. 
The developed OptiChat framework performs automatic error correction at the modeling level. 
In the latest stage of research, external validation mechanisms are further integrated with RL. 
A novel RL-based reward mechanism was introduced in AutoDH~\cite{ma2024automatic}, which balances solution quality, time cost and computational overhead of LLM invocation, enabling adaptive selection between expert-crafted and LLM-produced heuristic algorithms. 
Taking external optimization solvers as reward signal sources, Chen et al.~\cite{chen2025solver} put forward Solver-Informed RL to steer LLMs toward iterative optimization of problem-solving strategies according to validation feedback throughout the reinforcement learning workflow. 
Rather than deriving optimal values by solving optimization tasks, model equivalence is adopted to validate the precision of optimization modeling~\cite{wang2024optibench}. 
Additionally, Huang et al.~\cite{huang2025graphthought} generated effective inference sequences for solving graph combinatorial optimization problems using heuristic-guided forward search or backward reasoning strategies aligned with solvers.

\subsection[Internal Reflection]{Internal Reflection}

\begin{figure}[t]
	\centering
	\includegraphics[scale=0.30]{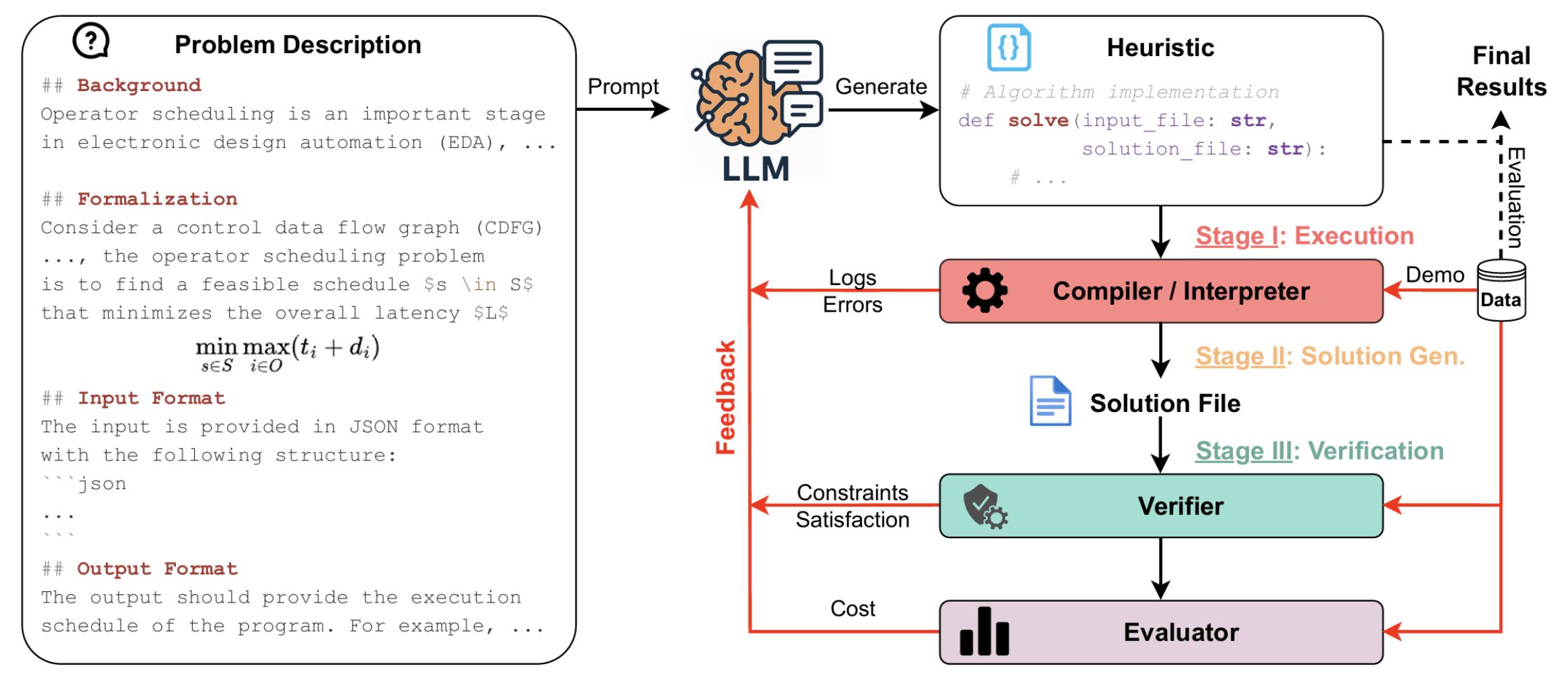} 
	\caption{Workflow of HeuriGym \cite{chen2026heurigym}.}

	\label{HeuriGym}
\end{figure}

Internal reflection methods focus on driving LLMs to identify potential errors and automatically execute corrective processes. 
A self-reflection and debugging mechanism was considered in~\cite{huang2024words}, allowing LLMs to autonomously detect flawed reasoning logic and perform corresponding revisions. 
Focusing on OR scenarios, Zhang et al.~\cite{zhang2025or} put forward an intelligent agent named OR-LLM-Agent, which realizes a closed cycle covering code running and error adjustment in sandbox settings. 
Integrating cluster analysis and search space introspection, Fornies et al.~\cite{fornies2025remoh} developed a multi-objective heuristic reflection-based evolutionary approach termed REMoH, which directs LLMs to yield varied and superior heuristic rules. 
In evolutionary computing research, Van et al.~\cite{van2025code} verified that LLMs stimulated by multi-turn prompts are capable of producing more satisfactory solutions. 
Equipped with majority voting and self-corrective reasoning tactics, OptiMind~\cite{OptiMind}  was established to facilitate LLMs in outputting more promising results. 
Apart from that, Chen et al.~\cite{chen2026heurigym} presented a platform called HeuriGym, which builds an iterative loop covering code generation, environmental execution and error feedback. Its operating process is illustrated in Fig.~\ref{HeuriGym}.

\subsection[Agent Collaboration]{Agent Collaboration}

In multi-agent systems, verification tasks are handled by dedicated agents, serving as independent components for automatic solution evaluation and code feasibility checking. 
For example, Elhenawy et al.~\cite{elhenawy2024visual} validated the effectiveness of multi-modal, multi-agent systems in combinatorial optimization tasks within a visual reasoning scenario. 
By conducting knowledge fusion and algorithmic solving for hierarchical modeling, Yuan et al.~\cite{yuan2025ma} put forward a method termed MA-GTS, which gradually reconstructs graph structures from textual content and invokes optimization algorithms in an adaptive manner. 
To explore agent-based heuristic optimization, Yang et al.~\cite{yang2025heuragenix} established a solution paradigm named HeurAgenix, where intelligent agents are capable of assessing newly generated and evolved heuristic strategies. 
Drawing on human cognitive patterns, an OR oriented system called ORMind~\cite{wang2025ormind} was constructed, supporting multi-round introspection over LLM-derived solutions to pursue optimal outcomes. 
Intelligent agents are adopted to pinpoint defects in logical reasoning and deliver revised prompts, whereby LLMs can polish their reasoning routes accordingly~\cite{tang2025calm}. 
Lima et al.~\cite{lima2025toward} proposed a mechanism that realizes cross verification via multi-agent voting. \\

\textit{In summary, existing solution validation frameworks  include solver-assisted external validation, self-consistency internal reflection, and multi-agent collaboration through division of labor, collectively forming the LLM4OR validation system.}

 \begin{figure}[t]
	\centering
	\includegraphics[scale=0.45]{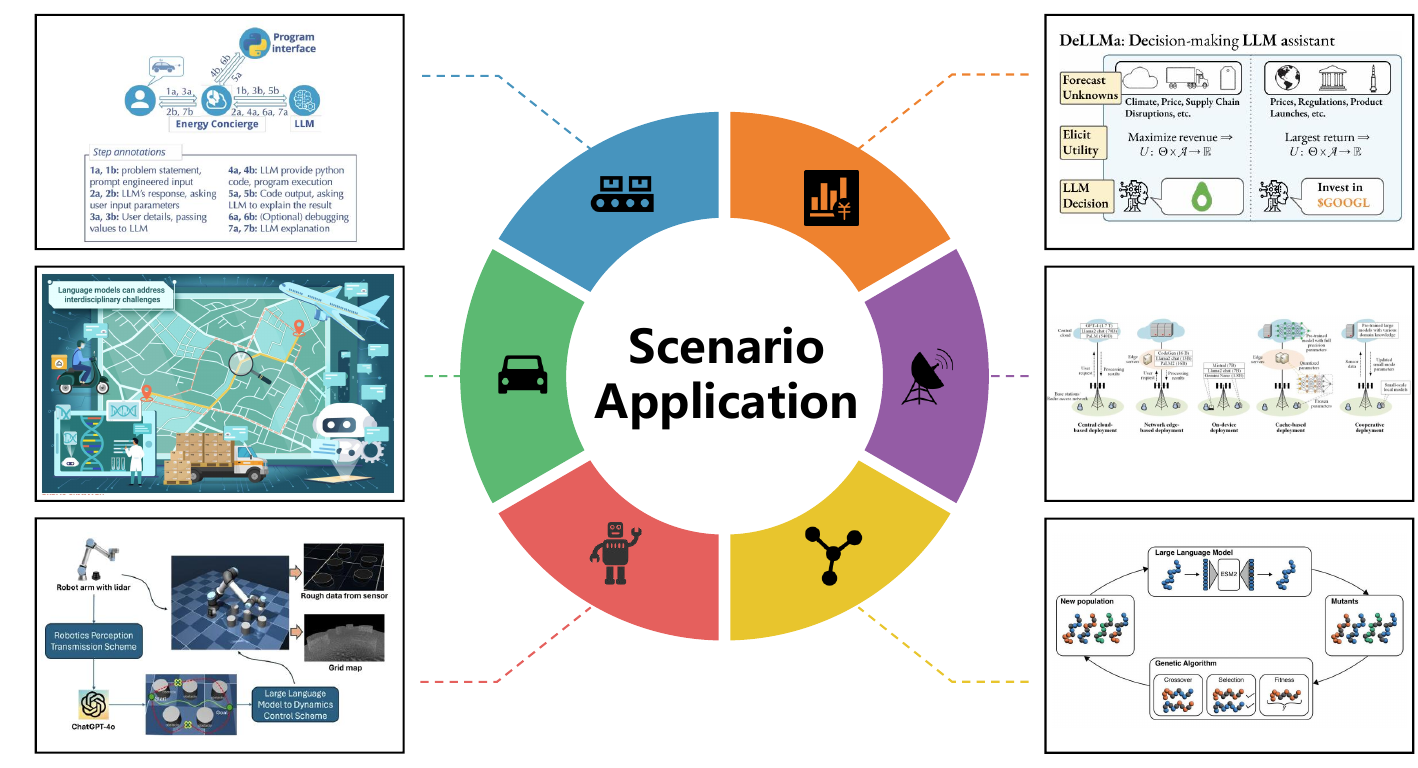} 
	\caption{Illustration of scenario application, such as scheduling~\cite{jin2024democratizing}, transportation~\cite{liu2023can}, robotics~\cite{mo2024precision}, biochemistry~\cite{nana2025integrating}, networks~\cite{zhou2024large}, and finance~\cite{liu2024dellma}}
	\label{ScenarioApplication}.
\end{figure}

\section[Scenario Application]{Scenario Application}\label{sec6}

This section reviews the representative applications of LLM4OR across different task categories, including scheduling, transportation, robotics, biochemistry, networks, finance, and other emerging domains, see Fig.~\ref{ScenarioApplication}.

\subsection[Scheduling]{Scheduling}

Scheduling is one of the most typical applications of OR, involving resource allocation and task sequencing. 
In recent years, LLMs have been used to generate scheduling strategies and automate solution finding under complex constraints~\cite{ivanov2026agentic}. 
To apply LLMs in supply chain scenarios, Li et al.~\cite{li2023large} designed OptiGuide, which was validated in a real Microsoft cloud supply chain case. 
Lawless et al.~\cite{lawless2024want} used LLMs for meeting scheduling to facilitate interactive decision support. 
Related studies have also explored the application of LLMs in program scheduling~\cite{jobson2024investigating} and task scheduling~\cite{yatong2024ts}. 
Combining RAG and CoT, Li et al.~\cite{li2026a4ps} devised a novel approach capable of adjusting tasks, revising constraints and upgrading algorithms for advanced planning and scheduling problems in manufacturing scenarios. 
Other research focuses on using LLMs to solve shop-floor scheduling problems. 
Jin et al.~\cite{jin2024democratizing} applied LLMs in energy network management, including tasks such as electric vehicle charging scheduling, HVAC control, and renewable energy planning. 
Abgaryan et al.~\cite{abgaryan2024llms} introduced a supervised dataset of 120,000 samples specifically for training LLMs to solve job shop scheduling problems, and fine-tuned the Phi-3-Mini model using LoRA.
To enhance the ability of LLMs in designing heuristic scheduling rules, Huang et al.~\cite{huang2024automatic} proposed  a population self-evolution strategy called SeEvo, which can automatically design heuristic algorithms for solving dynamic job shop scheduling problems.

\subsection[Transportation]{Transportation}

Transportation is a key application area for LLMs in intelligent mobility and transport. 
Capable of comprehending multi-source information, LLMs can cooperate with traditional optimization solvers to build decision-making systems that convert natural language demands into executable schemes. 
Some studies focus on using LLMs to improve the solving efficiency of the TSP and Vehicle Routing Problems (VRP)~\cite{liu2023can, tran2025large}. 
Other research combines methods such as MCTS, RL, or FT to mitigate issues like spatial hallucinations that LLMs may encounter in long-term planning, achieving performance improvements~\cite{zheng2025monte, ju2024globe}. 
In addition, the academic community has started exploring the integration of map data collection or external tool invocation into frameworks to build end-to-end decision systems that consider multi-objective constraints~\cite{da2024open}. 
By adopting fine-tuned LLMs to translate natural language requirements into symbolic forms, Ju et al.~\cite{ju2024globe} developed  TTG by integrating mixed-integer linear programming to derive optimal paths and completing the whole procedure within several seconds. 
To expand the external tool invocation capacity of existing models, Da et al.~\cite{da2024open} put forward a scheme named Open-TI, enabling end-to-end traffic decision-making ranging from map data collection to practical execution control. 
To efficiently solve TSP and VRP, Tran et al.~\cite{tran2025large} proposed a LLM-guided attention bias method to enhance the generalization capability of neural combinatorial optimization models.

\subsection[Robotics]{Robotics}

In the field of robotics, the introduction of LLMs enables collaborative planning and control, offering the advantage of combining environmental descriptions for resource allocation and task management for single or multiple robots. 
You et al.~\cite{you2023robot} applied LLMs in construction robot assembly to achieve automatic planning and dynamic adjustment of construction sequences. 
Dhanaraj et al.~\cite{dhanaraj2024preference} explored strategies to integrate human preferences into task planning by converting preferences into soft constraints, enabling human-robot collaboration. 
Obata et al.~\cite{obata2024lip}  combined linear programming with dependency graph structures, which achieved good results in experiments solving multi-robot task coordination and allocation. 
To explore the feasibility of LLMs in robot trajectory planning, Mo et al.~\cite{mo2024precision} used GPT-4o to provide real-time path planning algorithms, showing that their method outperforms existing methods in real-time feedback and trajectory accuracy. 
Huang et al.~\cite{huang2024can} adopted LLMs to realize self-debugging and self-verification, and achieved remarkable performance in single and multi-robot path planning tasks. 
In accordance with user instructions, Ji et al.~\cite{ji2026genswarm} proposed a strategy paradigm named GenSwarm, which is able to autonomously generate and deploy control policies adaptable to practical multi-robot systems.

\subsection[Biochemistry]{Biochemistry}

LLMs have become an important tool in facilitating biomolecule design, showing strong capabilities in sequence understanding and structure prediction. They are well-suited for integration with evolutionary algorithms or genetic search to efficiently explore biomolecular sequence spaces. 
Reinhart et al.~\cite{reinhart2024large} used LLMs as evolutionary optimizers to design large molecules with defined sequences. 
By applying LLMs to directed evolution processes, Tran et al.~\cite{tran2025protein} proposed MLDE, which outperforms existing benchmark algorithms in terms of efficiency and effectiveness in conditional protein generation tasks. 
Combining LLMs with genetic algorithms, Nana et al.~\cite{nana2025integrating} developed the enzyme design strategy, which first trains LLMs to learn the interrelationships between key amino acid residues in proteins, and then uses genetic algorithms to search for improved sequences, enhancing enzyme catalytic efficiency.

\subsection[Networks]{Networks}

Optimization problems in networks and communication often involve complex resource allocation, and reasoning capabilities of LLMs bring new automated solution approaches to this field~\cite{zhou2024large}. 
Targeting wireless network resource allocation scenarios, Qiu et al.~\cite{qiu2024large} presented a combinatorial optimizer termed LMCO, which seeks the optimal layout and access point deployment for wireless networks. 
For wireless communication resource scheduling, LLMs were utilized  to automatically identify nonconvex components and convert them into solvable expressions, thus realizing fully automated nonconvex resource allocation~\cite{peng2025llm}.

\subsection[Finance]{Finance}

LLMs have injected new vitality into the field of financial investment, offering better adaptability to new portfolio models and the dynamic market environment compared to traditional methods. 
To evaluate the effectiveness of LLMs in financial decision-making tasks, Liu et al.~\cite{liu2024dellma} introduced DeLLMa, which obtains good performance in real financial data experiments. 
Additionally, Li et al.~\cite{li2024generative} proposed a multi-objective evolutionary algorithm strategy based on LLMs, and validated superior performance in experiments compared to traditional multi-objective evolutionary algorithms. 
Yu et al.~\cite{yu2025finmem} constructed an LLM-driven multi-agent mode for financial trading, which achieves promising results in core indicators including cumulative return and Sharpe ratio. 
By integrating environmental and social factors into financial decision-making, Roy et al.~\cite{roy2026humanized} established a human-machine interactive intelligent analysis mode to boost decision-making and risk assessment capacity.

\subsection[Other Applications]{Other Applications}

Beyond the aforementioned fields, LLMs have also been extensively adopted in diverse general OR tasks. For example, to ensure that text generated by LLMs meets both logical and structural requirements, Regin et al.~\cite{regin2024combining} explored the possibility of embedding LLMs into constraint programming. 
Leveraging the logical reasoning capabilities of LLMs and integrating end-to-end models, Jiao et al.~\cite{jiao2024city} developed City-LEO to enhance the efficiency and transparency of smart city management. 
By combining network topology that characterizes engineering system component relationships with domain knowledge, Jiang et al.~\cite{jiang2025large} proposed an LLM-based framework to address engineering combinatorial optimization problems. In addition, Thomas et al.~\cite{thomas2025what} designed a universal mechanism integrating LLMs with combinatorial optimization methods, which significantly reduces the carbon emissions in menus while ensuring customer meal satisfaction.

\section[Datasets]{Datasets}\label{sec7}

To evaluate the capability of LLMs in solving OR problems, the research community has continuously expanded datasets and evaluation frameworks across various problem types and application domains. 
This section summarizes these relevant datasets, classifying them into two broad categories: general-purpose and task-specific.

\begin{table}[t]
\centering
\caption{General-purpose benchmark datasets.}
\label{tab2}
\renewcommand{\arraystretch}{1.2}
\footnotesize
\setlength{\tabcolsep}{4pt}
\begin{tabular*}{\textwidth}{@{\extracolsep\fill} l c l c}
\toprule
Datasets & Scales & Problem Types & Year \\
\midrule
NL4Opt \cite{ramamonjison2022augmenting} & 245 & LP & 2022 \\
NLP4LP \cite{ahmaditeshnizi2024optimus0.2} & 67 & LP, MILP & 2024 \\
NL2OPT \cite{mostajabdaveh2024optimization} & 70 & LP, MILP, QP & 2024 \\
ComplexOR \cite{xiao2024chain} & 37 & MILP & 2024 \\
MAMO \cite{huang2025llms} & 1,209 & LP, MILP, ODE (Easy \& Complex) & 2025 \\
IndustryOR \cite{huang2025orlm} & 100 & LP, IP, MILP, NLP & 2025 \\
EOR \cite{zhang2025decision} & 30 & LP, IP, MILP, NLP & 2025 \\
OptiBench \cite{yang2025optibench} & 605 & LP, IP, MILP, NLP & 2025 \\
DP-BENCH \cite{zhou2025auto} & 132 & DP & 2025 \\
CP-BENCH \cite{michailidis2025cp} & 101 & CP & 2025 \\
Text2Zinc \cite{singirikonda2025text2zinc} & 100 & LP, MILP, CP & 2025 \\
Bench4Opt \cite{wang2025orgeval} & 394 & LP, MILP & 2025 \\
OptMATH-Bench \cite{lu2025optmath} & 165 & LP, IP, MILP, NLP, SOCP & 2025 \\
\bottomrule
\end{tabular*}
\vspace{2pt}
\begin{flushleft}
\footnotesize
\textit{Note:} LP = Linear Programming, IP = Integer Programming, MILP = Mixed Integer Linear Programming, 
NLP = Nonlinear Programming, QP = Quadratic Programming, SOCP = Second-Order Cone Programming, 
CP = Constraint Programming, DP = Dynamic Programming, ODE = Ordinary Differential Equations.
\end{flushleft}
\end{table}

\subsection[General-Purpose Benchmark Datasets]{General-Purpose Benchmark Datasets}

General-purpose benchmark datasets focus on classical OR models, covering fundamental problem types such as linear programming, integer programming, and mixed-integer linear programming, as shown in Table~\ref{tab2}. 

Ramamonjison et al.~\cite{ramamonjison2022augmenting} released NL4Opt and established a standard evaluation paradigm in this field. 
To address the limitations of NL4Opt in modeling complex problems, Ahmaditeshnizi et al.~\cite{ahmaditeshnizi2024optimus0.2} introduced the NLP4LP dataset, which significantly increases problem complexity by incorporating long-text descriptions. 
Subsequently, evaluation frameworks have been extended to more complex scenarios. 
Xiao et al.~\cite{xiao2024chain} proposed the ComplexOR dataset by integrating academic papers and real-world industrial cases.

\begin{table}[t]
\centering
\caption{Accuracy comparison of representative studies on selected datasets.}
\label{tab3}
\renewcommand{\arraystretch}{1.2}
\footnotesize
\setlength{\tabcolsep}{4pt}
\begin{tabular*}{\textwidth}{@{\extracolsep\fill} l c c c c c}
\toprule
Methods
& \makecell{NL4Opt} 
& \makecell{IndustryOR} 
& \makecell{ComplexOR} 
& \makecell{MAMO\\(Easy)} 
& \makecell{MAMO\\(Complex)} \\
\midrule
CoE-GPT3.5-turbo \cite{xiao2024chain} & 58.9\% & - & 25.9\% & - & - \\
OptiMUS 0.2-GPT4 \cite{ahmaditeshnizi2024optimus0.2} & 78.8\% & - & 66.7\% & - & - \\
LEAN-LLM-OPT-GPT4.1 \cite{liang2026llm} & 94.7\% & 65.0\% & - & 93.5\% & 71.4\% \\
ORLM-Qwen2.5-7B \cite{huang2025orlm} & 86.1\% & 25.0\% & - & 85.2\% & 44.1\% \\
LLMOPT-Qwen1.5-14B \cite{jiang2025llmopt} & 93.0\% & 46.0\% & - & 97.0\% & 68.0\% \\
OptMATH-Qwen2.5-7B \cite{lu2025optmath} & 94.7\% & - & - & 86.5\% & 51.2\% \\
OR-R1-Qwen3-8B \cite{ding2026or} & 88.3\% & 35.3\% & 46.3\% & 86.1\% & 49.9\% \\
\bottomrule
\end{tabular*}
\end{table}

Dataset construction has gradually expanded to cover more problem types and become closer to real-world applications. 
Mostajabdaveh et al.~\cite{mostajabdaveh2024optimization} developed the NL2OPT dataset for practical applications. 
Huang et al.~\cite{huang2025llms} proposed the MAMO dataset, which evaluates solution quality by invoking solvers and comparing with ground truth, covering 1,209 problems including ordinary differential equations, linear programming, and mixed-integer linear programming. 
Huang et al.~\cite{huang2025orlm} further introduced the IndustryOR dataset, specifically designed to evaluate LLMs in solving real-world OR problems. 
Focusing on applications such as supply chain, financial investment, and logistics management, Zhang et al.~\cite{zhang2025decision} developed the EOR dataset based on IndustryOR.  To address the limitation that existing datasets lack tabular data, Yang et al.~\cite{yang2025optibench} proposed the OptiBench dataset. 
Some studies focus on specific problem types. 
Zhou et al.~\cite{zhou2025auto} introduced the DP-BENCH dataset for dynamic programming, while Michailidis et al.~\cite{michailidis2025cp} developed the CP-BENCH dataset for constraint programming. 
Singirikonda et al.~\cite{singirikonda2025text2zinc} designed the Text2Zinc dataset, which integrates satisfiability and optimization problems, enabling unified representation across domains. 
Additionally, Wang et al.~\cite{wang2025orgeval} proposed the graph-theoretic ORGEval benchmark, transforming model equivalence checking into a graph isomorphism problem, and constructed the accompanying Bench4Opt dataset, thereby enhancing the validation of modeling correctness. 
Lu et al.~\cite{lu2025optmath} introduced second-order cone programming problems in the OptMATH-Bench dataset, enabling the evaluation of LLMs in modeling more complex convex optimization problems.

As shown in Table~\ref{tab3}, when problems become more realistic or complex, performance generally declines significantly. 
Therefore, how to effectively improve the solution capability of existing methods remains a key focus for future research.

\subsection[Task-Specific Benchmark Datasets]{Task-Specific Benchmark Datasets}

Task-specific benchmark datasets focus more on complex OR problems, often involving tasks such as path planning and production scheduling. 
Table~\ref{tab4} presents the collected task-specific datasets. 

Some datasets focus on path planning problems. 
Wang et al.~\cite{wang2023can} proposed the NLGraph dataset, which contains 29,370 problems covering eight types of graph optimization tasks with varying complexity, such as topological sorting and Hamiltonian paths, establishing a foundation for evaluating the generalization capability of LLMs across multiple OR tasks. 
To assess the spatial and temporal reasoning capabilities of LLMs, Aghzal et al.~\cite{aghzal2024can} introduced the path planning dataset called PPNL. 
By constructing the RoutBench  dataset with 1,000 vehicle routing problems, Li et al.~\cite{li2025ars} designed the ARS framework to systematically evaluate the robustness of LLMs in generating heuristic algorithms. 
Based on the AtCoder Heuristic Contest, Imajuku et al.~\cite{imajuku2025ale} released the ALE-BENCH dataset, which covers tasks such as path planning and production scheduling.

\begin{table}[t]
\centering
\caption{Task-specific benchmark datasets.}
\label{tab4}
\renewcommand{\arraystretch}{1.2}
\footnotesize
\setlength{\tabcolsep}{4pt}
\begin{tabular*}{\textwidth}{@{\extracolsep\fill} l p{9.5cm} l}
\toprule
Datasets & Domains & Year \\
\midrule
NLGraph \cite{wang2023can} & 29,370 graph-structured optimization problems covering 8 task types & 2023 \\
PPNL \cite{aghzal2024can} & Path planning tasks & 2024 \\
asyncHow \cite{lin2024graphenhanced} & Problem instances from domains such as diet, education, and healthcare & 2024 \\
RoutBench \cite{li2025ars} & Including 1,000 VRP instances and their variants & 2025 \\
ALE-BENCH \cite{imajuku2025ale} & NP-hard problems such as path planning and production scheduling & 2025 \\
GraphArena \cite{tang2025grapharena} & 4 polynomial-time challenges and 6 NP-hard optimization tasks & 2025 \\
CO-BENCH \cite{sun2025co} & 36 types of real-world combinatorial optimization tasks & 2025 \\
FrontierCO \cite{feng2025comprehensive} & 8 representative combinatorial optimization tasks & 2025 \\
OPT-BENCH \cite{li2025opt} & 20 real-world ML tasks and 10 NP-hard combinatorial optimization problems & 2025 \\
HeuriGym \cite{chen2026heurigym} & 9 domains including chip design, protein engineering, and electronic circuits & 2026 \\
\bottomrule
\end{tabular*}
\end{table}

An increasing number of studies evaluate LLM performance across broader applications. 
Lin et al.~\cite{lin2024graphenhanced} released the asyncHow dataset, which includes 1,600 high-quality instances from domains such as diet, education, and healthcare. 
Tang et al.~\cite{tang2025grapharena} constructed the GraphArena dataset, containing 10,000 instances, including four types of polynomial-time challenges such as common neighbors and six NP-hard optimization tasks such as the TSP. 
Sun et al.~\cite{sun2025co} introduced the CO-BENCH dataset, which covers 36 types of complex optimization tasks such as packing, cutting, and layout, and establishes a standardized pipeline spanning data processing, algorithm design, and performance evaluation. 
In addition, Feng et al.~\cite{feng2025comprehensive} expanded evaluation dimensions by covering eight representative combinatorial optimization tasks in the FrontierCO dataset. 
To simulate more complete scientific and engineering problem solving processes, Li et al.~\cite{li2025opt} constructed the OPT-BENCH dataset, which includes five categories of ML subtasks such as regression and prediction, as well as typical NP-hard problems in graph theory, set selection, and scheduling. 
Chen et al.~\cite{chen2026heurigym} constructed the HeuriGym dataset covering nine domains such as chip design, protein engineering, and airline operations.

\section[Conclusion]{Conclusion}\label{sec8}

The development of LLMs has provided new perspectives for solving OR problems. 
This paper systematically reviews the latest research progress of LLMs in the field of OR, with a particular focus on their key roles throughout the entire problem solving and application pipeline. 
To end this survey, we analyze several challenging issues  and outline potential directions for future research.

\subsection[Framework Exploration]{Framework Exploration}

As shown in Table~\ref{tab3}, the accuracy on some complex datasets is below 50\%. 
Moreover, existing studies are mostly restricted to small-scale datasets or local performance improvements, making it difficult to handle industrial-scale problems with millions of variables, external CSV parameters, and multi-stage stochastic structures. 
Undoubtedly, framework exploration targeting complex decision-making scenarios has become a key research trend in LLM4OR~\cite{xiao2026deepor, xie2026murka}. 
In addition, developing domain-specific LLMs for particular OR applications is also necessary.

\subsection[Dataset Construction]{Dataset Construction}

Existing benchmark datasets are primarily derived from textbooks, which makes it difficult to reflect the complexity and diversity of real-world OR scenarios. In addition, although various benchmark datasets have been proposed to evaluate LLM capabilities, they lack unified design standards, making comprehensive and consistent comparisons challenging. 
Accordingly, effective datasets construction and unified evaluation criteria have become essential measures to advance the overall performance of LLM4OR in practical research~\cite{chen2026chain}.

\subsection[Data Protection]{Data Protection}

Industrial data are often highly sensitive, such as road network information in transportation logistics and transaction records in financial investment. This poses significant challenges for cloud-based LLMs in model formulation and application deployment. 
Given the risk of privacy leakage during data transmission, reasonable data protection has become a key prerequisite to explore offline deployment schemes and privacy-preserving strategies, greatly promoting the practical popularization of LLM4OR technologies~\cite{zhang2026collaborative, veera2026securing}.

\subsection[Result Explainability]{Result Explainability}

The interpretability of LLMs not only facilitates experts in debugging generated results but also provides insights for improving performance. 
For example, Explainable Operations Research (EOR)~\cite{zhang2025decision} proposed evaluation methods to analyze causal relationships between modeling decisions and optimization results, while the state-chain semantic alignment framework~\cite{wang2025translating} integrated LLM-based generation methods with symbolic theorem provers. 
In this regard, result explainability is a core yet challenging research priority for LLM4OR.

\subsection[Multimodal LLMs]{Multimodal LLMs}

OR problems often involve not only textual descriptions but also multiple data modalities such as images and tables. 
However, solutions relying on LLMs still exhibit limitations in handling heterogeneous data, thereby constraining their performance. 
This has motivated researchers to explore the potential advantages of Multimodal Large Language Models (MLLMs) in problem solving~\cite{chen2025mathflow, zhao2026underappreciated}. 
Designing application frameworks based on MLLMs is expected to provide more efficient solutions for decision-making in complex environments~\cite{wang2026graphcogent}.

\section*{Acknowledgements}
This work was supported by the National Natural Science Foundation of
China under Grant 12371306 and 12571345.





\bibliography{mybibfile}

\end{document}